\documentclass[11pt,a4paper]{article}
\usepackage{amsmath,amssymb,pifont,enumerate,latexsym}

\long\def\comment#1{} 
\def \qed {{\hspace*{\fill}\ding{111}\\[2mm]}}

\def \card {\mathrm{card}\,} 
\def \dom {\mathrm{dom}\,}

\def \T  {\mathcal{T}} 
\def \complex    {{\mathrm{Complex}\,}}

 \newcommand{\s}{\rightarrow}

\def \bbbn{{\mathrm{I\!N}}} 
\def \N  {\bbbn}

\newenvironment{proof}{\textit{Proof. }}{\qed} 
 
\newtheorem{theorem}{Theorem}[section]

\newtheorem{corollary}[theorem]{Corollary} 
 
\newtheorem{lemma}[theorem]{Lemma}

\newtheorem{proposition}[theorem]{Proposition} 

\begin{document} 
\title{\bf Is Complexity a Source of Incompleteness?} 
\author{{\bf Cristian S. Calude}\\
Department of Computer Science\\
The University of Auckland\\ Private Bag 92019,
Auckland, New Zealand\\ 
Email: {\tt cristian@cs.auckland.ac.nz} \\ 
\and {\bf  Helmut J{\"u}rgensen}\\
Department of Computer Science,
 The University of Western Ontario\\ London,
Ontario, Canada N6A 5B7\\
  Institut f\"{u}r Informatik, Universit\"{a}t
Potsdam, August-Babel-Str. 89\\   D-14469 Potsdam, Germany\\
} \date{}
 
\maketitle 
\noindent 
\thispagestyle{empty} 
 
\begin{abstract} 
In this paper we prove   Chaitin's ``heuristic principle'',  
{\it the theorems of a finitely-specified  theory cannot be significantly more complex than the theory itself},  for an appropriate measure of complexity.  We show that the measure is invariant under the change of the G\"odel numbering.  For this measure, the theorems of a 
 finitely-specified, sound, consistent theory  strong enough to formalize arithmetic  which is arithmetically sound  (like Zermelo-Fraenkel set theory with choice  or  Peano Arithmetic)  have bounded complexity, hence every  sentence of the theory  which is significantly more complex than the theory is unprovable.  Previous results showing that  incompleteness is not accidental, but ubiquitous are here reinforced in probabilistic terms:  
  the probability that a true sentence of length  $n$ is provable in the theory  tends to zero when $n$ tends to infinity, while the  probability that   a sentence of length  $n$  is true  is strictly positive.  
\end{abstract}

\section{Introduction}
G\"odel's Incompleteness Theorem states that every  finitely-specified, sound, theory which is strong enough to include arithmetic cannot be both consistent and complete. G\"odel's original proof as well as most  subsequent proofs are based on the following idea:  a theory which is finitely-specified, sound, consistent and  strong enough can express sentences about provability within the theory,  which,  themselves,  are not provable by the theory, but  can be shown to be true using  a proof by contradiction.
A true and unprovable sentence is called independent.
 This type of proof of incompleteness  does not answer the questions
of whether independence is a widespread phenomenon nor which kinds of sentences can be expected to be independent.

 Chaitin \cite{ch3}   presented a complexity-theoretic proof of G\" odel's Incompleteness Theorem  which shows that high complexity is a reason of the unprovability of  infinitely many  (true)  sentences. This result
 suggested to him  the following  ``heuristic principle'', a kind of information-preservation principle: 
{\it  the theorems of a finitely specified  theory cannot be significantly more complex than the theory itself}.
This approach would address the second of the questions above, that is, highly complex sentences are  independent, and, as a consequence, would  indicate that independence is pervasive.
A formal confirmation of the pervasiveness of independence has been obtained in \cite{cjz} via a topological analysis; a quantitative result is still missing.

In this paper we prove that a formal version of  the ``heuristic principle'' is indeed correct for an appropriate measure of complexity; the  measure  is invariant under the change of the G\"odel numbering. 
For this measure, $\delta$, the theorems of a  finitely-specified, sound, consistent  theory which is strong enough to include arithmetic have bounded complexity, hence every sentence of the theory  which is significantly more complex than the theory is unprovable.  Previous results showing that  incompleteness is not accidental, but ubiquitous are here reinforced in probabilistic terms:  
  the probability that a true sentence of length  $n$ is provable in the theory  tends to zero when $n$ tends to infinity, while the  probability that   a sentence of length  $n$  is true  is strictly positive.

The paper is organized as follows. In Sections~2 and~3 we present the background,  the notation and main results needed for our proofs. In Section~4  we  discuss some general complexity-theoretic results which will be used to prove  the main  result (Theorem~\ref{complex-incompleteness}). In Section~5 we prove that incompleteness  is  widespread in probabilistic terms. In Section~6  we use the new complexity measure to  prove Chaitin's information-theoretic incompleteness result for the Omega Number.
We finish with a few general comments  in Section~7. The bibliography includes a selection of relevant papers and books, but  is by no means complete.

\section{Background}

G\"odel's Incompleteness Theorem, announced on 7 October 1930  in K\"onigsberg at 
 the First International Conference on the Philosophy of Mathematics\footnote{Hilbert, von Neumann, Carnap, Heyting presented reports; the conference was part of the
 German Mathematical Congress.}  is a landmark of the  twentieth century mathematics (see \cite{godel,heijenoort,feferman2} for the original paper, \cite{nn,Smor,cp,hintikka} for other proofs and more related mathematical facts, 
\cite{kleene,kreisel,rodrig,dawson,wang,def,RozenbergSalomaa,ca2,castijim,hintikka,crismind} for more mathematical, historical and philosophical details). It says that {\it  in a   finitely-specified, sound, consistent  theory  strong enough to formalize arithmetic, there are true, but unprovable sentences;}  so such a theory is {\it  incomplete}. A true and unprovable sentence is called {\it independent}.   The first condition states that  axioms can be algorithmically listed;  consistency means free of contradictions;   soundness means that any proved  sentence is true.

According to Hintikka (\cite{hintikka}, p.~4),  with the exception of  von Neumann, who immediately grasped G\"odel's line of thought and its importance,    incompleteness passed un-noticed in  K\"onigsberg: even  the speaker who summarized the discussion  omitted G\"odel's result.
In spite of being praised, discussed, used or abused by many authors, the Incompleteness Theorem  seems, even after so many years since its discovery, stranger than most mathematical theorems.\footnote{This is quite impressive, as mathematics abounds with  baffling results.}
For example, according to Solovay (\cite{Ko},
p.~399):   ``The feeling
was that G\"odel's theorem was of interest only to logicians''; 
in Smory\'nski's words, (\cite{Ko}, p.~399),  ``It is fashionable to deride G\"odel's theorem
as artificial, as dependent on a linguistic trick.''

 In 1974 Chaitin \cite{ch3}   presented a complexity-theoretic proof of G\" odel's Incompleteness Theorem  which shows that high complexity is a reason of the unprovability of  infinitely many  (true)  sentences. This complexity-theoretic approach was discussed by Chaitin \cite{ch75,ch7,ch8,ch10,ch11,ch12} and
 various authors including  Davis \cite{davis},  Tymoczko \cite{ty},  Boolos and Jeffrey \cite{boolos},  pp.~288--291,  Svozil \cite{sv},  Li and  Vit\'{a}nyi \cite{lv},  Barrow \cite{barrow,barrow0},
 Calude \cite{ca1,Ca}, Calude and Salomaa \cite{caludesalomaa}, Casti \cite{casti},
Delahaye  \cite{jeanpaul}; it was critized by  
 van Lambalgen \cite{la2}, Fallis  \cite{fallis}, Raatikainen \cite{raat}, Hintikka \cite{hintikka}.
 
  Chaitin's proof in  \cite{ch3} is based on program-size  complexity (Chaitin complexity) $H$:   the complexity $H(s)$ of a binary string $s$  is the size, in bits, of the shortest program for a universal self-delimiting Turing machine  to calculate $s$. The  complexity $H(s)$  is unbounded. The proof shows that 
 {\it for every finitely-specified, sound, consistent theory  strong enough to formalize arithmetic, there exists
 a positive constant $M$ such that no  sentence of the form ``$H(x)>m$'' is provable in the theory unless $m$ is less than $M$.}  There are infinitely many true  sentences of the form  ``$H(x)>m$'' with
 $m >M$, and each of them  is unprovable in the theory.

 The high $H$-complexity of the sentences ``$H(x)>m$'' with
 $m >M$  is a source of their unprovability.\footnote{Fallis \cite{fallis}, p.~264,  argued that G\"odel's true but unprovable sentence G  is likely to have excessive $H$-complexity. Similarly, if the theory is capable of expressing its own consistency, then the corresponding sentence is likely to have excessive $H$-complexity. It would be interesting to have a mathematical confirmation of  these facts.}
 Is every true sentence $s$ with $H(s)>M$ unprovable by the theory? Unfortunately, the answer is {\it negative}  because only finitely many sentences $s$ have complexity $H(s) <M$
 in contrast with the fact that the set of all  theorems of the theory is infinite. For example, 
  $ZFC$ (Zermelo-Fraenkel set theory with choice)  or  Peano Arithmetic trivially
 prove all  sentences of the form ``$n+1=1+n$''. The $H$-complexity
 of the sentences  ``$n+1=1+n$'' grows unbounded with $n$. This fact, noticed and discussed by Chaitin in Section 6 of \cite{ch13}  (reprinted in \cite{ch11}  pp.~55--81) as well as by  Svozil  \cite{sv},  pp.~123--125,   is  essential for  the critique  in   \cite{fallis,raat} (cited in \cite{hintikka}); the works  \cite{ch13,ch11,sv} seem to be unknown to the  authors of \cite{fallis,raat,hintikka}.

Chaitin's proof based on $H$  cannot be directly extended to all  unprovable sentences, hence the problem of whether complexity is a source of incompleteness remained  open. In this note we prove that  the ``heuristic principle''  proposed by Chaitin in  \cite{ch11}, p.~69,  namely that
{\it the theorems of a  finitely-specified  theory cannot be significantly more complex than the theory itself} \footnote{An  ``approximation'' of this principle supported by Chaitin's proof is that
``one cannot prove, from a set of axioms, a theorem that is of greater $H$-complexity than the axioms {\it and know} that one has done it''; see \cite{ch11}, p.~69; see also Theorem~\ref{chi} in Section~4.} is correct if  we measure the complexity of a string by the difference between the program-size  complexity and the length of the string, our $\delta$-complexity (Theorem~\ref{complex-incompleteness}).  The $H$-complexity
 of the sentences  ``$n+1=1+n$'' grows unbounded with $n$, but the ``intuitive complexity'' of the sentences  ``$n+1=1+n$'' remains  bounded; this intuition is confirmed by  $\delta$-complexity. Note that a sentence with a large $\delta$-complexity has also a large $H$-complexity, but the converse is not true. There are only {\it finitely} many strings with bounded $H$-complexity, but {\it infinitely} many strings
 with bounded $\delta$-complexity.
 
As a consequence of  Theorem~\ref{complex-incompleteness} we  prove that   the incompleteness phenomenon is more widespread than   previously shown in \cite{godel,feferman2,ch3,ch10,ch11} and by the  topological analysis of  \cite{cjz}:
   the probability that a true sentence of length  $n$ is provable in the theory  tends to zero when $n$ tends to infinity, while the  probability that   a sentence of length  $n$  is true  is strictly positive.

\section{Prerequisites}\label{s.prel} 

We  follow the notation in \cite{Ca}. By 
$\N=\{0,1,2,\ldots \}$ we denote the set of non-negative integers.
Further on, $ \log_{Q}$ denotes the base $Q\ge 2$ logarithm and
   $ \log n = \lfloor  \log_{2} (n+1) \rfloor $; $ \lfloor \alpha \rfloor$ is 
the ``floor" of the real $\alpha$ and  $\lceil \alpha \rceil$ 
is the ``ceiling" of  $\alpha$.
The cardinality
of the set $A$ is denoted by $\card(A)$.  An alphabet  with $Q$ elements will be denoted by $X_{Q}$; by
$X_{Q}^*$ we denote the set of finite strings (words) on $X_{Q}$, including
the {\itshape empty\/} string $\lambda$.  The length of the string $w\in X_{Q}^*$
is denoted by $|w|_{Q}$. 

For $Q=2$ we use the special  set  $B = \{0,1\}$ instead of $X_{2}$.
We consider the following bijection between non-negative integers and strings on  $B$:
$ 0 \,\, \mapsto \,\, \lambda,
1 \,\, \mapsto \,\, 0,
 2 \,\, \mapsto \,\, 1, 3 \,\, \mapsto \,\, 00, 4 \,\, \mapsto \,\, 01,  5 \,\, \mapsto \,\, 10,  6 \,\, \mapsto \,\, 11,  \ldots $
The image of $n$, denoted $bin(n)$, is the binary representation
of the number $n+1$ without the leading 1. Its length is  $| bin(n) |_{2} = \log n$.
In general we denote by
$string_Q (n)$
 the $n$th string on  $X_{Q}$  according to the quasi-lexicogra\-phi\-cal order.
In particular, $bin(n) = string_{2}(n)$. In this 
way we get a bijective function  $string_Q: \N  \s X_{Q}^{*}$;
$|string_{Q}(n)|_{Q} = \lfloor \log_{Q} (n(Q-1)+1) \rfloor$.

We assume that the reader is familiar  with  Turing machines processing
 strings, computability and program-size  complexity (see, for example,
\cite{ca1,boolos,Ca,downeybook}). The program set (domain) of the Turing machine $T$ is the set
$PROG_T= \{x\in X_{Q}^*  :  T \mbox{ halts on } x \}$; when $T$ halts on $x$, $T(x)$ is the result of the computation of $T$ on $x$.
A partial function $\varphi$ from strings to strings is called 
partial computable  (abbreviated p.c.) if there is a Turing machine $T$ such that:   a)  $PROG_T = \dom(\varphi)$, and b) $T(x)=\varphi(x)$, for each
$x\in PROG_T$.  
A computable function is a p.c.\ function  $\varphi$ with $ \dom(\varphi) = X_{Q}^{*}$.
A set of strings is computable if its characteristic function is computable.  A set of strings is  computably 
 enumerable (abbreviated c.e.)  if it is  the program set of a  Turing machine. 

A {\it self-delimiting Turing machine}  is a Turing machine $T$ such that its program set  is a prefix-free set of strings.  Recall that a {\em prefix-free set} of strings $S$  is a set such that no string in $S$  is a proper extension of any other string in $S$.  In what follows the term  {\it machine}  will refer to either  a p.c.\ function with prefix-free domain or a  self-delimiting Turing machine.

 Each prefix-free set  $S\subset X_{Q}^{*}$ satisfies  Kraft's inequality:   $\sum_{i=1}^{\infty} r_{i}\cdot Q^{-i} \le 1$,  where $r_{i} = \card \{x\in S  :  |x|_{Q}=i\}$.  A stronger result,
 the Kraft-Chaitin Theorem (see \cite{Ca}, p.~53),  is essential in algorithmic information theory:   Let $n_{1}, n_{2}, \ldots $ be a computable sequence of non-negative integers such that
\begin{equation}
\label{kraft}
\sum_{i=1}^{\infty} Q^{-n_{i}} \le 1.
\end{equation}
 Then, we can effectively construct a prefix-free sequence of strings $w_{1}, w_{2}, \ldots $ such that for each $i\ge 1, \, |w_{i}|_{Q}=n_{i}$.

The {\it program-size  complexity\/} of the string $x\in X_{Q}^*$ (relative   to $T$)
is $H_{Q,T}(x)=\min \big\{|y|_{Q}  :  y \in X_{Q}^*, \ T(y)=x \big \}$, where 
$\min \emptyset = \infty$.  The Invariance Theorem states that we can effectively construct a  machine $U=U_{Q}$ (called {\it universal$\,$}) such that for every
machine $T$ there exists a constant $\varepsilon >0$ such that for all $x\in X_{Q}^{*}$,   $H_{Q,U} (x) \leq H_{Q,T} (x) + \varepsilon $. In what follows we will
fix $U$ and put $H_{Q} = H_{Q,U}$; in particular, $H_{2}$ denotes the program-size  complexity induced by a universal (binary)  machine. If $x$ is in $X^{*}_{Q}$, then $x^{*} = \min
\{u\in X^{*}_{Q}: U_{Q}(u) = x\}$, where the minimum is taken according to the quasi-lexicographical order; it is seen that $H_{Q}(x) = |x^{*}|_{Q}$.

\section{Complexity and Incompleteness}

In this section we introduce the   $\delta$-measure and then prove for it  Chaitin's ``heuristic principle'': 
{\it  the theorems of a  finitely-specified  theory cannot be significantly more complex than the theory itself}.

First we introduce the $\delta$-measure.
Recall that $U_{Q}$ is a fixed universal machine on $X_{Q}$ and  $H_{Q} = H_{Q,U_{Q}}$.
In what follows we will work with the function $\delta_{Q} (x) =  H_{Q}(x)-|x|_{Q}$ (note that $-\delta_{Q}$ is a   ``deficiency of randomness'' function in the sense of  \cite{Ca}, Definition 5.21, p.~113). The  $\delta$-complexity  is ``close'', but not equal,  to the conditional $H_{Q}$-complexity, of a string given its length.

The complexity measures $H_{Q}$ and $\delta_{Q}$ have  similarities as $\delta_{Q}$ is defined from $H_{Q}$ by means of some simple computable functions;  for example,  they are both uncomputable. But $H_{Q}$ and $\delta_{Q}$ differ in an {\it essential}  way: given a positive $N$, the set
$\{x\in X_{Q}^{*}  :   H_{Q}(x) \le N\}$ is finite while, by Corrolary~\ref{infinity},  the set $\{x\in X_{Q}^{*}  : \delta_{Q}(x) \le  N\}$ is infinite.  Note that the conditional $H_{Q}$-complexity does not have this property.  A sentence with a large $\delta_{Q}$-complexity has also a large $H_{Q}$-complexity, but the converse is not true. For example,  the $H_{Q}$-complexity of a (true) sentence of the form  ``$1+n=n+1$'' is  about  $\lfloor \log_{Q} n \rfloor$ plus a constant, a function which tends to infinity as $n \s \infty$;  however,  their $\delta_{Q}$-complexity is bounded.

In view of Theorem~5.4 in \cite{Ca}, p. 102, there exists a constant $c>0$ such that:

\begin{equation}
\label{deltamax}
\max_{|x|_{Q}=N} \delta_{Q}(x) \ge  H_{Q}(string_{Q}(N)) -c,
\end{equation}

so there are strings of arbitrarily large $\delta_{Q}$-complexity.

\bigskip

 The following result  is taken from \cite{Ca} (Theorem 5.31, p.~117). 

\begin{theorem}
\label{immune}
 For every $t \geq 0$, the set  $C_{Q,t} = \{x \in  X_{Q}^{*} :  \delta_{Q}(x) >-t\}$  is immune, that is,  the set  is infinite and contains no infinite c.e.\ subset.
\end{theorem}
\if01
\begin{proof}
Clearly, the set  $C_{Q,t}$ is infinite as,   for each $N$, $\card \{ x \in  X^{*}_{Q}  :  |x|=N, \delta(x)>-t\} $ 
$>
2^{N} (1- 2^{-t}) \ge 0$. Next, assume  that  $C_{Q,t}$  contains an infinite c.e.\ subset, hence  it contains an infinite computable subset  $D$. For $x,y\in X^{*}_{Q}$, we put $x < y$ if $bin^{-1}(x)<bin^{-1}(y)$, and define the machine $T$     by 
$T(0^{i}1) =  \min\{ x \in D  :   |x| \geq t+2(i+1) \}$. Clearly, $T$ has a  computable graph. Because $T(0^{i}1)\in D\subset C_{Q,t}$,   for each $i \ge 1$,  we have:
$$H(T(0^{i}1))  > |T(0^{i}1)|-t
                            \geq t+2(i+1)-t
                            =2(i+1).$$
Consequently, by the Invariance Theorem,  there exists a constant $\varepsilon >0$ such that  for  infinitely many   $i$, we have
$$2(i+1) < H(T(0^{i}1)) \leq H_{T}(T(0^{i}1))+\varepsilon \leq i+1+\varepsilon,$$  a contradiction.
\end{proof}
\fi

\begin{corollary}
\label{complex}
 For every $t \geq 0$, the set  $\complex_{Q,t} = \{x \in  X^{*}_{Q}  :   \delta_{Q}(x) >t\}$  is immune.
\end{corollary}
\begin{proof} 
As $\complex_{Q,t} \subset C_{Q,t}$ and every  infinite subset of an immune set is immune itself, we only need to show that $\complex_{Q,t}$ is infinite. To this aim we use  formula  (\ref{deltamax}) and the fact that the function
 $H_{Q}(string_{Q}(N))$ is unbounded.
\end{proof}

\begin{corollary}
\label{infinity}
For every $t \ge 0$, the set 
$\{ x\in X_{Q}^{*} : \delta_{Q}(x) \le t\}$  is infinite.
\end{corollary} 
\begin{proof} 
 The set in the statement is not even  c.e.\  because, by Corollary~\ref{complex}, its complement  is immune.
 \end{proof}

The above result suggests that any ``reasonable'' theory cannot include more than finitely many theorems with high $\delta$-complexity. And, indeed, a simple analysis confirms this fact. A formal language used by a theory capable of speaking about natural numbers includes
variables (a fixed variable  $x$ and the sign $^{\prime}$ may be used to generate all variables,   $x, x^{\prime}, ^{\prime\prime}$, etc.), the constant 0, function symbols for successor, addition and multiplication,
$s, +, \cdot,$  the sign for equality, =, logical connectives,$\neg, \wedge, \vee, \Rightarrow$, quantifiers, $\forall, \exists$, and parentheses, (,). They form an  alphabet $X_{15}$.\footnote{Of course, we can work with smaller or larger alphabets, depending on specific needs.} The formal language consists of well-formed formulae which respect strict syntactical rules;
for example, each left parenthesis has to be matched with exactly one right parenthesis.
Theorems are then defined by specifying the axioms and the inference rules. For instance, the  system {\bf Q} introduced by R. M. Robinson (see, for example,
\cite{ec}) contains the logical axioms (propositional, substitution, $\forall$-distribution, equality axioms) and the following seven axioms: Q1: $(s(x) = s(x^{\prime}))
\Rightarrow  (x = x^{\prime}),$ Q2: $ \neg (0 =s(x))$, Q3: $(\neg (x=0)) \Rightarrow
\exists x^{\prime} (x =s(x^{\prime}))$, Q4: $x+0=x$, Q5: $x +s(x^{\prime}) = s(x+x^{\prime})$, Q6: $x\cdot 0 =0$, Q7: $x\cdot s(x^{\prime}) = (x\cdot x^{\prime}) + x$,
and the inference rules of  modus ponens and generalisation. A proof in the system {\bf Q} is a sequence of well-formed formulae such that each formula is either an axiom, or is derived from two earlier  formulae in the sequence by an inference rule.    {\it Theorems} are well-formed formulae which have proofs in {\bf Q}.
As theorems are special well-formed formulae, it is clear that each theorem $x$
in the system {\bf Q} has rather small $H_{15}$-complexity, more precisely, $H_{15}(x)$ is not larger than its length plus a fixed constant. Such a remark suggests that
Chaitin's ``heuristic principle'' may be true for $\delta_{15}$. However, this property could be a consequence of some particular way of writing/coding the theorems!  
To be able to measure somehow
the ``intrinsic'' complexity of a theorem   we  need to prove that the property 
is  {\it  invariant with respect to a system of acceptable names}, in our case,
{\it G\"odel numberings}.

To make the discussion precise, let us fix a formal language $L\subset X_{Q}^{*}$. A {\it G\"odel numbering} for $L$ is a computable, one-to-one function $ g: L \s B^{*}$, i.e. a system of unique binary names for the well-formed formulae of $L$. For example, a G\"odel numbering for the well-formed formulae of the system {\bf Q} can be obtained by coding
the elements of the alphabet $X_{15}$ with the first 15 binary strings of length four, and then extend this coding according to the syntax of the language. Various other possibilities can be imagined;  see for example, \cite{boolos,ec}.

As the set of theorems is a c.e. subset of the set of well-formed formulae, we will work only with computable, one-to-one functions $ g: \T \s B^{*}$
defined on the set of theorems.

{\it The $\delta$-complexity of a theorem $u\in\T$ induced by the G\"odel numbering} $g$
is defined by: 
\begin{equation}
\label{deltadef}
\delta_{g}(u) = H_{2}(g(u)) -\lceil \log_{2}Q\rceil  \cdot |u|_{Q}.
\end{equation}

 The  formula for  $\delta_{f}$   is essentially  the formula defining $\delta_{Q}$  relativized to the  G\"odel numbering $g$: the factor  $\lceil \log_{2}Q\rceil  $ has the role of ``adjusting'' the sizes of the alphabets $X_{Q}$ and $B$.

 The first result confirms the  intuition: we prove that $\delta_{g}$  is, up to an additive constant,
  equal to  $\lceil \log_{2}Q\rceil  \cdot \delta_{Q}$.
 
 \begin{theorem}
 \label{invat}
Let $\T\subset X_{Q}^{*}$ be c.e.\ and $g : \T  \s B^{*}$ be  a G\"odel numbering. Then, there effectively exists a constant $c$ (depending upon  $U_{Q}, U_{2}$, and $g$) such that for all $u\in \T$ we have:
\begin{equation}
\label{inva}
|\delta_{g}(u) - \lceil \log_{2}Q\rceil  \cdot \delta_{Q}(u) |\le  c.
\end{equation}
\end{theorem}
\begin{proof}
First we prove the existence of a constant $c_{1}$ such that 
\begin{equation}
\label{inva1}
\delta_{g}(u) \le \lceil \log_{2}Q\rceil  \cdot \delta_{Q}(u) +  c_{1}.
\end{equation}

For each string $w \in PROG_{U_{Q}}$ we define
$n_{w} =  \lceil \log_{2}Q\rceil \cdot |w|_{Q},$
and we note that
\[\sum_{w \in  PROG_{U_{Q}}} 2^{-n_{w}} = \sum_{w \in  PROG_{U_{Q}}} 2^{- \lceil \log_{2}Q\rceil \cdot |w|_{Q}} \le \sum_{w \in  PROG_{U_{Q}}} Q^{-|w|_{Q}} \le 1,\]
because $ PROG_{U_{Q}}$ is prefix-free. Using now the Kraft-Chaitin Theorem, we can effectively construct, for every $w\in PROG_{U_{Q}}$ a binary string $s_{w}$ such that
$|s_{w}|_{2} = n_{w}$ and the set $\{s_{w} : w \in  PROG_{U_{Q}}\}$ is c.e.\ and prefix-free. This allows us to construct the  machine $C$ defined by
\[C(s_{w}) = g(U_{Q}(w)), \mbox{  for  } w\in PROG_{U_{Q}}.\]
As $C(s_{w^{*}}) = g(U_{Q}(w^{*}))=g(w)$ we have $$H_{C}(g(w)) \le |s_{w^{*}}|_{2} =
 \lceil \log_{2}Q\rceil \cdot |w^{*}|_{Q} =  \lceil \log_{2}Q\rceil \cdot H_{Q}(w).$$ Applying the Invariance Theorem we get a constant $c_{1}>0$
such that
\begin{eqnarray*}
\delta_{g}(w) &=& H_{2}(g(w)) -  \lceil \log_{2}Q\rceil \cdot |w|_{Q}\\
& \le & 
 \lceil \log_{2}Q\rceil \cdot (H_{Q}(w) - |w|_{Q}) + c_{1} \\
 &=&  \lceil \log_{2}Q\rceil \cdot \delta_{Q}(w)+c_{1},
 \end{eqnarray*}
 which proves (\ref{inva1}).
 
 Secondly we prove the existence of a constant $c_{2}$ such that
 \begin{equation}
\label{inva2}
 \lceil \log_{2}Q\rceil  \cdot \delta_{Q}(u) \le \delta_{g}(u) +  c_{2}.
\end{equation}

For each $w\in PROG_{U_{2}}$ such that $|w|_{2} \ge \log_{2}Q$, we put
$m_{w}=  \lceil |w|_{2}\cdot \log_{Q}2 \rceil\ge 1,$
and we note that

\[\sum_{w\in PROG_{U_{2}}, |w|_{2} \ge \log_{2}Q} \, Q^{-m_{w}} \le \sum_{w\in PROG_{U_{2}}, |w|_{2} \ge \log_{2}Q}  \, 2^{-|w|_{2}} \le 1,\]

hence, in view of the   Kraft-Chaitin Theorem, we can effectively construct, for every $w\in PROG_{U_{2}}$ with $ |w|_{2} \ge \log_{2}Q$, a string $t_{w}\in X^{*}_{Q}$ of length $|t_{w}|_{Q}=m_{w}$ such that
the set $\{t_{w} : w\in PROG_{U_{2}}\}$ is c.e.\  and prefix-free.  In this way we construct the  machine $D$ defined by $D(t_{w}) = u$ if $U_{2}(w) = g(u)$. This construction is well-defined because $g$ is a G\"odel numbering. It is seen that if $U_{2}(w) = u$ and
$|w|_{2} \ge \log_{2}Q$, then $H_{D}(u) \le \lceil |w|_{2}\cdot  \log_{Q}2 \rceil$, so
applying the Invariance Theorem we get a constant $d$ such that 
\[\lceil \log_{2}Q \rceil \cdot
H_{Q}(u) \le   \lceil \log_{2}Q \rceil \cdot
H_{D}(u) + d \le H_{2}(g(u)) + d,\]  hence there is a constant $c_{2}$ such that
(\ref{inva2}) becomes true. We have used the fact that $\lceil \log_{2}Q \rceil \cdot
 \lceil m\cdot \log_{Q}2 \rceil \le m$, for all integers  $m>0$.

Finally, (\ref{inva}) follows  from (\ref{inva1}) and (\ref{inva2}).
\end{proof}

 As a consequence,
 asymptotically,  the $\delta$-measure is independent of the G\"odel numbering.

\begin{corollary}
Let $\T\subset X_{Q}^{*}$ be c.e. and $g, g^{\prime} : \T \s B^{*}$ be two  G\"odel numberings.
 Then, there effectively exists a constant $c$ (depending upon
 $ U_{2}, g$ and $ g^{\prime}$) such that for all $u\in \T$ we have:
\begin{equation}
\label{invg}
|\delta_{g}(u) - \delta_{g^{\prime}}(u)| \le c.
\end{equation}
\end{corollary}
\begin{proof}
The relation  (\ref{invg}) follows from Theorem~\ref{invat}. However, it is instructive to give a short, direct proof. To this aim consider the  machine $C$ defined for $w\in B^{*}$ by $C(w) = g(u)$ if $U_{2}(w) = g^{\prime}(u)$.
The definition is correct because $PROG_{C} \subset PROG_{U_{2}}$ and $g$ is computable and one-to-one. If $U_{2}(s) = g^{\prime}(u)$, then $C(s) = g(u)$, so by the 
Invariance Theorem there exists a constant $c_{1}$ such that for all $u\in L$,
$\delta_{g}(u) \le \delta_{g^{\prime}}(u) + c_{1}.$ Finally, (\ref{invg}) follows by symmetry.
\end{proof}

\begin{theorem}
\label{complex-incompleteness}
Consider   a  finitely-specified, arithmetically sound (i.e. each arithmetical proven sentence is true),  consistent theory  strong enough to formalize arithmetic,
and denote by $\T$ its set of theorems written in the alphabet $X_{Q}$.  Let $g$ be a G\"odel numbering for $\T$. Then, there  exists 
a constant $N$,  which depends upon $U_{Q}, U_{2}$ and $ \T$,  such that 
 $\T$  contains no  $x$  with $\delta_{g}(x)>N$. 
\end{theorem}
\begin{proof}
 Because of syntactical constraints, there exists a positive constant $d$ such that for every $x\in\T$, $H_{Q} (x)
\le |x|_{Q} +d$, i.e. $\delta_{Q}(x) \le d$ (see also the  discussion of the system $\mathbf  Q$ following Corrolary~\ref{infinity}). Hence in view of Theorem~\ref{invat}, there is a constant $N\ge d$ such that
for every  $x\in\T$, $\delta_{g}(x) \le N$.
\end{proof}

Every sentence $x$ in the language of  $\T$  with  $\delta_{g}(x)>N$  is unprovable in the theory; every such ``true'' sentence is thus independent of the theory. 

Do we have examples of such sentences? First, Chaitin's sentences
of the form ``$H_{2}(x) > n$'', for large $n$ are  such examples. 

Here is another
way to construct true sentences of high $\delta$-complexity.
A formula $\varphi (x)$ in the language of arithmetic is called $\Sigma_{1}$ if it is of the form $(\exists y)\theta(x,y)$, where $\theta$  contains only two free variables $x$ and $y$.  We  write
$\N \models \varphi (n)$ to mean that $\varphi (n)$ is true when $n$ is interpreted as a non-negative integer.
The Representation Theorem (see \cite{Smor})  states that a set  $R\subset \N$ is c.e. iff there (effectively)  exists a $\Sigma_{1}$ formula $\varphi (x)$ such that for all $n \in \N$ we have:
$n \in R \, \Leftrightarrow \, \N \models \varphi (n).$

For every $a\in\N$, the set $\{n\in\N  :   \delta_{Q}(string_{Q}(n)) \le a\}$ is c.e., so in view of the Representation Theorem  there exists a $\Sigma_{1}$ formula $\varphi$ (depending on $U_{Q}, a$) such that for
every $n\in\N$ we have: $ \delta_{Q}(string_{Q}(n)) \le a \, \Leftrightarrow \, \N \models \varphi (n)$.
Consequently, the formula $\psi = \neg \varphi$  represents the predicate ``$\delta_{Q}(string_{Q}(n))>a$''. Because of consistency and soundness, by enumerating the theorems in $\T$ of the form $\psi(m)$ (corresponding to  true formulae $\psi(m)$) we get an enumeration of the set
$ \{x \in \T : \psi(string_{Q}^{-1}(x))\in \T\}  \subset \{x\in\T : \delta_{Q}(x)>a\}$.

Now let $a$ be a non-negative integer.  As  $ \{x \in \T : \psi(string_{Q}^{-1}(x))\in \T\}$ is a c.e. subset of the immune set $\{x\in X^{*}_{Q}  :  \delta_{Q}(x)>a\}$,  it has to be finite, that is, there exists an $M\in\N$ such that for every $x\in\T$ with $\psi(string_{Q}^{-1}(x)) \in\T$ we have:  $|x|_{Q}\le M$.  
We have got Chaitin's statement (\cite{ch11}, p.~69): 

\begin{theorem}
\label{chi}
Every finitely-specified,  arithmetically sound,  consistent theory  strong enough to formalize arithmetic  can prove only, for finitely many of its  theorems, that they have high $\delta$-complexity.
 \end{theorem}
 
 The theory can formalise all sentences of the form $\psi(m)$ in a very economical way, i.e. with small $\delta$-complexity, but is incapable of proving more than finitely many instances: almost all
 true formulae of the form $\psi(m)$ are unprovable.

{\bf Comments\/} (a)  
  Theorem~\ref{complex-incompleteness} establishes a limit on the $\delta_{g}$-complexity of provable sentences in $\T$;    the bound depends upon the chosen G\"odel numbering $g$. In this approach it makes no sense to measure the power of the theory by its complexity, i.e.  through the minimal $N$ such that the theory proves no  sentence $x$  with $\delta_{g}(x)>N$  (see also the discussion in \cite{la2}).

(b)  Theorem~\ref{complex-incompleteness}  does not hold true for an arbitrary finitely-specified theory
as there are c.e.\ sets containing strings of arbitrary large $\delta$-complexity.

(c)  It is possible to have incomplete theories without high $\delta$-complexity sentences; for example, an incomplete theory for propositional tautologies.

\section{Is Incompleteness  Widespread?}

The first application complements the result of  \cite{cjz} stating that the set of  unprovable sentences is topologically  large.  We probabilistically  show  that only a  few true sentences can be proven in a given theory, but   the set  of true sentences is ``large''.  

We begin with the following result:

\begin{proposition}
Let $N>0$ be a fixed integer,  $T\subset X^{*}_{Q}$ be c.e. and $g: T \s B^{*}$
be a  G\"odel numbering. Then,
\begin{equation}
\label{l}
\lim_{n \rightarrow \infty} Q^{-n} \cdot \card{\{x\in X^{*}_{Q} : |x|_{Q}=n, \delta_{g}(x) \le N\}}=0.
\end{equation}
\label{limzero}
\end{proposition}
\begin{proof} We present here a direct proof.\footnote{Alternatively, one can  evaluate
 the size of the set of strings of a given length having almost maximum $\delta_{Q}$-complexity.} In view of Theorem~\ref{invat}, there exists a constant $c>0$ such that
 \[ \{x\in X^{*}_{Q} : |x|_{Q}=n, \delta_{g}(x) \le N\} \subseteq
  \{x\in X^{*}_{Q} : |x|_{Q}=n,  \lceil \log_{2}Q \rceil \cdot  \delta_{Q}(x) \le N + c\}.\]
 So, we only need to evaluate the limit\\[-5ex]
 
 \begin{equation}
\label{ll}
\lim_{n \rightarrow \infty} Q^{-n} \cdot \card{\{x\in X^{*}_{Q} : |x|_{Q}=n, \delta_{Q}(x) \le M\}}=0,
\end{equation}

 where $M = \lceil (N+c)/ \lceil \log_{2}Q \rceil \rceil$.
 
 First we note that for every $n$ we have: $
 \{x\in X^{*}_{Q} : |x|_{Q}=n, \delta_{Q}(x) \le M\}
=  \{x\in X^{*}_{Q} : |x|_{Q}=n,  \exists y\in X_{Q}^{*} \,  (|y|_{Q} \le n+M, U_{Q}(y) = x)\}$, 
so\\[-5ex]

\begin{eqnarray*}
\card{  \{x\in X^{*}_{Q} : |x|_{Q}=n, \delta_{Q}(x) \le M\}}
 & \le & \card{ \{y\in X^{*}_{Q}  : |y|_{Q} \le n + M, |U_{Q}(y)|_{Q}=n\}}\\
 & \le & \card{ \{y\in X^{*}_{Q}  : |y|_{Q} \le n + M, U_{Q}(y) \mbox{  halts}\}}\\
 \end{eqnarray*}\\[-5ex]
Consequently, \\[-5ex]

\begin{equation}
\label{limit}
\lim_{n \rightarrow \infty} Q^{-n} \cdot \card{  \{x\in X^{*}_{Q} : |x|_{Q}=n, \delta_{Q}(x) \le M\}} =
\lim_{n \rightarrow \infty} \sum_{i=1}^{n+M} Q^{-n} \cdot r_{i},
\end{equation}

where $r_{i} = \card \{y\in X^{*}_{Q} : |y|_{Q}=i, U_{Q}(y) \mbox{ halts}\}.$ Using the Stolz-Cesaro Theorem we get
\begin{equation}
\label{limcal}
\lim_{n \rightarrow \infty} \sum_{i=1}^{n+M} Q^{-n} \cdot r_{i} = Q^{M} \cdot \lim_{m \rightarrow \infty} \sum_{i=1}^{m} Q^{-i} \cdot r_{i} =
 Q^{M}/(Q-1) \cdot \lim_{m \rightarrow \infty} Q^{-m}\cdot r_{m} =0,
 \end{equation} 
  due to Kraft's inequality  $ \sum_{i=1}^{\infty} r_{i}\cdot Q^{-i} \le 1$. So, in view of , (\ref{ll}), (\ref{limit}) and (\ref{limcal}) we get (\ref{l}).
\end{proof}

\begin{theorem}
\label{prob}
Consider   a consistent, sound,  finitely-specified theory  strong enough to formalize arithmetic. 
The probability that a true  sentence of length  $n$ is provable in the theory  tends to zero when $n$ tends to infinity,  while the  probability that   a sentence of length  $n$  is true  is strictly positive.  
\end{theorem}

\begin{proof}
We fix    a consistent, sound,  finitely-specified theory, let  $\T$ be its set of theorems
and let $g$ be a G\"odel numbering of $\T$. For every integer $n\ge 1$, let $\T^{n} =
\{x\in \T : |x|_{Q}=n\}$. By Theorem~\ref{complex-incompleteness}, there exists a positive integer $N$ such that   $\T  \subseteq \{x\in X^{*}_{Q} :   \delta_{g}(x) \le N\}$. Consequently, for every $n$:  
$ \T^{n}  \subseteq  \{x\in X^{*}_{Q} :  |x|_{Q}=n, \delta_{g}(x) \le N\},$ so in view of Proposition~\ref{limzero}, the  probability that a sentence of length  $n$ is provable in the theory  tends to zero when $n$ tends to infinity.

Next consider the sentences $h_{x,m}= $``$H_{Q}(x)>m$'', where $x$ is a string over the alphabet $X_{Q}$. For every $n\ge 1$ there exists a positive integer $N_{m}$ such that
for every string $x\in X_{Q}^{*}$ of length $|x|_{Q}> N_{m}$, $h_{x,m}$ is true. 

For each fixed $m$, $|h_{x,m}|_{Q} = |x|_{Q}+c$, so for every $m\ge 1$ and $n \ge N_{m}+c$ we have:
\[
\card\{w\in X^{*}_{Q}: |w|_{Q} = n, w \mbox{  is true}\} \cdot Q^{-n} \ge 
\card\{x\in X^{*}_{Q}: |x|_{Q} = n-c, \} \cdot Q^{-n} \ge  Q^{-c},
\]
showing that the  probability that  a sentence of length $n$  is true  is strictly positive. 
\end{proof}

\section{Incompleteness and $\Omega_{U}$}

The second application is to use $\delta$ to prove Chaitin's  Incompleteness Theorem for $\Omega_{U}$ \cite{ch75} (see also the analysis in \cite{jeanpaul,crissolovay,Ca}).  This shows that $\delta$ is a ``natural'' complexity.
We start with the following preliminary result: 

\begin{lemma}
\label{limitcomplex}
 Let $x_{1}x_{2}\cdots $ be an infinite binary  sequence and let $F$ be a strictly increasing function mapping positive integers to positive integers. If the set $\{(F(i), x_{F(i)})  :  i \ge 1\}$ is c.e., then there exists a constant $\varepsilon >0$ (which depends upon $U$ and the characteristic function of the  set) such that for all $k\ge 1$ we have:
\begin{equation}
\label{bc}
\delta_{2}(x_{1}x_{2}\cdots x_{F(k)}) \le \varepsilon -k.
\end{equation}
\end{lemma}
\begin{proof}
To prove (\ref{bc}), for $k\ge 1$ we consider the strings:
\begin{equation}
\label{s}
w_{1}x_{F(1)}w_{2}x_{F(2)}\cdots w_{k}x_{F(k)},
\end{equation}
where each $w_{j}$ is a string of length $F(j)-F(j-1)-1, F(0)=0$. In this way
we effectively generate all binary strings of length $F(k)$ in which    the bits
on the ``marked''positions $F(1), \ldots ,F(k)$ are fixed.

It is clear that $\sum_{i=1}^{k} |w_{i}| = F(k)-k$ and the mapping $(w_{1}, w_{2},  \ldots ,w_{k}) \mapsto w_{1}w_{2} \cdots w_{k}$ is bijective, hence to generate all strings of the form (\ref{s}) we only need to generate all strings of length $F(k)-k$. Hence, we
 consider the enumeration of all strings of the form (\ref{s}) for $k=1,2,\ldots$ The lengths of these strings  form the sequence
\[ \underbrace{F(1), F(1), \ldots ,F(1)}_{2^{F(1)-1}\mbox{ times}}, \ldots
\underbrace{,F(k), F(k), \ldots ,F(k)}_{2^{F(k)-k}\mbox{  times}}, \ldots\]
which is computable and satisfies the inequality (\ref{kraft}) as
$\sum_{k=1}^{\infty}  2^{F(k)-k}\cdot 2^{-F(k)} = 1.$
Hence, by the   Kraft-Chaitin Theorem, for every string $w$ of length $F(k)-k$ there effectively exists a string $z_{w}$ having the same length as $w$ such that the set $\{z_{w} \in B^{*} : |w|_{2}=F(k)-k, k\ge 1\}$ is prefix-free. 
Indeed, from  a string $w$ of length $F(k)-k$  we get  a unique decomposition $w = w_{1}\cdots w_{k}$, and $z_{w}$ as above, so we can define   $C(z_{w}) = w_{1}x_{F(1)}w_{2}x_{F(2)}\cdots w_{k}x_{F(k)}$; $C$ is a  machine. Clearly,  $\delta_{C}(w_{1}x_{F(1)}w_{2}x_{F(2)}\cdots w_{k}x_{F(k)}) \le  |z_w|_{2} - F(k) = -k,$
for all $k\ge 1$. So by the Invariance Theorem we get the inequality (\ref{bc}).
\end{proof}

Consider now Chaitin's Omega Number, the halting probability of $U$:  $\Omega_{U} = 0.\omega_{1}\omega_{2}\cdots$, see \cite{ch4}. The binary sequence $\omega_{1}\omega_{2}\cdots$
is (algorithmically) random. There are various ways to characterize randomness (see for example \cite{ch8,Ca,downeybook}). A particular useful way is
the following complexity-theoretic criterion due to Chaitin: there exists a positive constant $\mu$ such that  for every $n\ge1$, 
\begin{equation}
\label{u}
\delta_{2} (\omega_{1}\omega_{2}\cdots \omega_{n}) \ge -\mu.
\end{equation}

The condition  (\ref{u}) is equivalent to $\sum_{n=0}^{\infty}  2^{- \delta_{2}  (\omega_{1}\omega_{2}\cdots \omega_{n})} < \infty$, cf, \cite{jm}.

It is easy to see that the inequality (\ref{bc})  in Lemma~\ref{limitcomplex} contradicts (\ref{u}), so a  sequence $x_{1}x_{2}\cdots x_{n}\cdots $ satisfying the hypothesis of  Lemma~\ref{limitcomplex}  cannot be  random. 

\begin{theorem}
Consider   a  consistent, sound,  finitely-specified theory  strong enough to formalize arithmetic. Then, we  can  effectively 
compute a constant $N$ such that the theory  cannot determine more than $N$ scattered digits of $\Omega_{U} = 0.\omega_{1}\omega_{2}\cdots$
\end{theorem}
\begin{proof}
Assume  that the theory  can determine infinitely many digits of $\Omega_{U} = 0.\omega_{1}\omega_{2}\cdots$
Then, we could effectively  enumerate an infinite sequence of digits of $\Omega_{U}$, thus satisfying the hypothesis of Lemma~\ref{limitcomplex} which would contradict the randomness of  $\omega_{1}\omega_{2}\cdots$
\end{proof}

\section{Conclusions}

  There are  various  illuminating proofs  of G\"odel's Incompleteness Theorem   and  some  interesting examples of true but unprovable sentences (see for example \cite{ch11,ph,Ha}).
   Still,  the phenomenon of incompleteness seems, even after almost 75 years since its  discovery, strange and to a large extent irrelevant to   `mainstream  mathematics', whatever this expression might mean.  Something is {\it missing}  from the picture. Of course, the  `grand examples' are missing; for example, no important open problem except Hilbert's tenth problem, see \cite{mat}, was proved to be unprovable. Other
 questions of interest include the source of incompleteness  and how common  the incompleteness phenomenon is. These two last questions  have been investigated in this note.

Chaitin's   complexity-theoretic proof of G\" odel's Incompleteness Theorem  \cite{ch3} shows that high complexity is a sufficient reason for the unprovability of  infinitely many  (true) sentences. This approach suggested that excessive complexity may  be a source of incompleteness, and, in fact, Chaitin (in \cite{ch13,ch11}) stated this as a ``heuristic principle'':   ``the theorems of a finitely-specified  theory cannot be significantly more complex than the theory itself''. By changing the measure of complexity, from program-size  $H(x)$ to $\delta(x)=H(x)-|x|$, we have proved (Theorem~\ref{complex-incompleteness}) that for any    finitely-specified, sound, consistent  theory  strong enough to formalize arithmetic  (like Zermelo-Fraenkel set theory with choice  or  Peano Arithmetic) and for any G\"odel numbering $g$ of its well-formed formulae,  we can compute a bound $N$ such that no sentence $x$ with complexity $\delta_{g}(x) >N$ can be proved in the theory; this phenomenon is independent on the choice of the G\"odel numbering.    For  a theory satisfying the hypotheses of Theorem~\ref{complex-incompleteness},    the probability that a true sentence of length  $n$ is provable in the theory  tends to zero when $n$ tends to infinity, while the  probability that   a sentence of length  $n$  is true  is strictly positive. This result reinforces the  analysis in \cite{cjz} which shows that the set of  independent  sentences is topogically  large.

According to Theorem~\ref{complex-incompleteness},  sentences expressed by  strings with  large $\delta_{g}$-complexity are unprovable. Is the converse implication true? In other words, given a theory as in the statement of Theorem~\ref{complex-incompleteness},
are there  independent
 sentences $x$  with low  $\delta_{g}$-complexity? Even if such sentences do exist,  in view of Theorem~\ref{prob}, the probability that
  a true sentence of length $n$ with $\delta_{g}$-complexity less than or equal to $N$ is 
  unprovable in the theory
  tends to zero when $n$ tends to infinity.

Other open questions which are interesting to  study include:  a)  the complexity of some concrete independent  sentences, like the  sentence expressing the consistency of the theory itself, b) the problem of finding other (more interesting?)  measures of  complexity satisfying Theorem~\ref{complex-incompleteness}, c) a stronger version of Theorem~\ref{prob}: under the same conditions, the probability that a sentence  of length  $n$, expressible in the language of the theory, is provable in the theory  tends to zero when $n$ tends to infinity.

\section*{Acknowledgement} We thank Henning Bordihn, John Casti, Greg Chaitin, Maia Hoeberechts and Andr\' e Nies for stimulating discussions and critical comments.

\end{document}